\documentstyle[12pt]{article}
\begin{document}
\begin{center}
{\Large \bf Characterization of Hermitian symmetric spaces }

\medskip
{\Large \bf by fundamental forms}

\bigskip
 {\large \bf Jun-Muk Hwang}\footnote{ Supported by the Korea Research Foundation Grant
(KRF-2002-070-C00003).}
 \hspace{3mm} {\large \bf and} \hspace{3mm} {\large \bf Keizo
 Yamaguchi}\footnote{Supported by Gran-in-Aid for Scientific Research (No.11304002),
Japan Society for the Promotion of Science.}
 \end{center}

\section{Introduction}

In [L1] and [L2], Landsberg proved the following

\medskip
{\bf Theorem A (Landsberg)} {\it Let $H$ be an irreducible compact
Hermitian symmetric space of rank 2, different from  the
hyperquadric ${\bf Q}_n \subset {\bf P}_{n+1}$. Let $H \subset
{\bf P}_N$ be a minimal non-degenerate equivariant embedding,
equivalently, an embedding of $H$ in ${\bf P}_N$ defined by the
complete linear system associated to the ample generator of
$Pic(H) \cong {\bf Z}$. Let $M \subset {\bf P}_N$ be a (not
necessarily closed) complex sub-manifold with $\dim(M) = \dim(H)$
and $x \in M$ be a point in a neighborhood of which all the
integer-valued differential invariants of $M$ remain constant. If
the second fundamental form of $M$ at $x$ is isomorphic to the
second fundamental form of $H$ at a point, then $M$ is
projective-linearly equivalent to an open subset of $H$.}

\medskip
For the reducible case, he proved ([L1, Theorem 1])

\medskip
{\bf Theorem B (Landsberg)} {\it Let $M \subset {\bf P}_N$ be a
complex sub-manifold and $x \in M$ be a point in a neighborhood of
which all the integer-valued differential invariants of $M$ remain
constant. If the second fundamental form of $M$ at $x$ is
isomorphic to the second fundamental form of a Segre variety ${\bf
P}_{n_1} \times {\bf P}_{n_2} \subset {\bf P}_{(n_1+1)(n_2+1)-1}$,
with $ n_1, n_2 \geq 2$, then $M$ is projective-linearly
equivalent to an open subset of the Segre variety.}

\medskip
We refer the readers to [GH] and [L1] for historical back-ground
on problems of this kind and  the relation of such questions with
other problems in projective algebraic geometry.

\medskip
Our main result is a generalization of Theorem A  to higher-order
embeddings and to Hermitian symmetric spaces of higher rank.
Precisely, we are going to prove

\medskip
{\bf Theorem 1} {\it Let $H$ be a  compact Hermitian symmetric
space which  contains neither ${\bf P}_n, n\geq 1$ nor $ {\bf
Q}_n, n \geq 2$ as an irreducible factor. Let $H \subset {\bf
P}_N$ be a non-degenerate equivariant embedding, or equivalently,
the embedding  by a complete linear system on $H$. Let $M \subset
{\bf P}_N$ be a (not necessarily closed) complex sub-manifold with
$\dim(M) = \dim(H)$ and $x \in M$ be a point in a neighborhood of
which all the integer-valued differential invariants of $M$ remain
constant. If the fundamental forms of $M$ at $x$ are isomorphic to
the fundamental forms of $H$ at a point, then $M$ is
projective-linearly equivalent to an open subset of $H$.}

\medskip
Note that for the minimal embedding of a Hermitian symmetric space
of rank 2, the second osculating space spans the whole ${\bf
P}_N$, so the higher fundamental forms are identically zero. Thus
Theorem A is a special case of Theorem 1.

 It is well-known that the equivariant embeddings of ${\bf
Q}_n \subset {\bf P}_{n+1}, n \geq 2$ and  ${\bf P}_n, n \geq 1$
cannot be characterized by their fundamental forms. For example,
for the standard embedding of ${\bf Q}_n$ in ${\bf P}_{n+1}$ and
the second Veronese embedding of ${\bf P}_n$, one has to consider
the cubic invariants in addition to the fundamental forms ([JM],
[L1, Theorem 4]). However, even when $H$ contains ${\bf P}_n$ or
${\bf Q}_n$ as an irreducible factor, there are  cases where our
method can be applied.  As an example, we will prove the following
generalization of Theorem B:

\medskip
{\bf Theorem 2} {\it Let $M \subset {\bf P}_N$ be a complex
sub-manifold and $x \in M$ be a point in a neighborhood of which
all the integer-valued differential invariants of $M$ remain
constant. If the  fundamental forms of $M$ at $x$ is isomorphic to
those of a Segre variety ${\bf P}_{n_1}  \times \cdots \times {\bf
P}_{n_m} \subset {\bf P}_N, N= \prod_{i=1}^m (n_i+1) \; -1, n_i
\geq 2$, then $M$ is projective-linearly equivalent to an open
subset of the Segre variety.}

\medskip
For the proof of Theorem 1 and Theorem 2, we use the theory of
linear differential equations modelled on Hermitian symmetric
spaces due to Y. Se-ashi [S]. This theory reduces the problem to
the vanishing of certain Spencer cohomology groups by a reduction
theory of principal bundles.  In this sense, the proof is
fundamentally different from [L1] and [L2], which depend on
certain moving frame computation. Of course, one may say that the
moving frame method and the reduction theory of principal bundles
are just two different formulations of the same procedure.
However, we believe that the more conceptual view-point of [S] is
essential  in the treatment of the higher rank  and the higher
order cases, because a simple-minded moving frame computation
would get complicated heavily.  We hope that presenting such a
geometric application would make the deep work [S] more widely
appreciated. See [SYY] for another application.

\medskip
In Section 2, we will define the fundamental forms  in the setting
of [S]. In Section 3, an invariant-theoretic result concerning the
fundamental forms of the Hermitian symmetric space is proved,
which makes it possible to apply the result of [S] in our problem.
This reduces the problem to a calculation of Lie algebra
cohomology.  The computation of cohomology group will be done
using Kostant's harmonic theory in Section 4. This computation had
already appeared in [SYY] in the irreducible case.

\section{Fundamental forms as a system of linear differential
equations}

We start with recalling some definitions from [S] (see also
[SYY]). Let $M$ be a connected complex manifold. We denote by $T$
and $T^*$ the (holomorphic) tangent and the cotangent bundle of
$M$ respectively. For a line bundle $E$ on $M$, we denote by
$J^p(E)$ the bundle of $p$-jets of $E$. Its fiber at a point $x
\in M$ is
$$J^p_x(E) := {\cal O}_x(E) / \left({\cal O}_x(E) \otimes {\bf
m}_x^{p+1} \right).$$ We identify $J^0(E) = E$.  For each pair of
non-negative integers $p \geq q,$ we have a natural projection
$\pi^p_q: J^p(E) \rightarrow J^q(E)$. For a section $s$ of $E$,
denote its $p$-jet at $x$  by $j^p_x(s) \in J^p_x(E)$. We have a
natural  exact sequence of vector bundles
$$ 0 \longrightarrow {\rm Sym}^p T^*(M) \otimes E \longrightarrow
J^p(E) \stackrel{\pi^p_{p-1}}{\longrightarrow} J^{p-1}(E)
\longrightarrow 0.$$ At each $x \in M$, we have the induced exact
sequence of vector spaces of the fibers at $x$ $$0 \longrightarrow
{\rm Sym}^p T^*_x \otimes E_x \longrightarrow J^p_x(E)
\stackrel{\pi^p_{p-1}}{\longrightarrow} J^{p-1}_x(E)
\longrightarrow 0.$$

A sub-bundle $R$ of $J^p(E)$ is called a {\bf system of linear
differential equations of order} $p$ on $E$. A  section $s$ of $E$
over an open subset $U \subset M$ satisfying $j^p_x(s) \in R_x$
for each $x \in U$ is called a {\bf solution} of $R$ over $U$. For
each non-negative integer $r \leq p$, the subspace $\sigma^r_x(R)
\subset {\rm Sym}^r T^*_x$ corresponding to $\pi^p_r(R_x) \cap
({\rm Sym}^r T^*_x \otimes E_x)$ under an identification $E_x
\cong {\bf C}$ is called the {\bf $r$-th  symbol} of $R$ at $x$.
It is independent of the choice of the identification $E_x \cong
{\bf C}$. The direct sum $\sigma_x(R) = \bigoplus_{r=0}^p
\sigma_x^r(R)$ is called the {\bf symbol} of $R$ at $x$. The
system $R \subset J^p(E)$ is said to be {\bf of finite type} if
$\sigma^p_x(R) =0$ for general $x \in X$. $R$ is said to be {\bf
integrable} if for each $\eta \in R_x$, there exists a  solution
$s$ over a neighborhood of $x$ for which $j^p_x(s) =\eta$. Such a
solution is uniquely determined by $\eta$ when $R$ is of finite
type.

Let $R$ (resp. $R'$) be a system of order $p$ on a line bundle $E$
(resp. $E'$) on a complex manifold $M$ (resp. $M'$). We say that
the symbol of $R$ at $x \in M$ is isomorphic to the symbol of $R'$
at $y \in M'$ if there exists  a linear isomorphism $\theta: T_x
\rightarrow T_y$ such that the induced linear isomorphism
$\theta^r: {\rm Sym}^r T^*_x \rightarrow {\rm Sym}^r T^*_y$ sends
$\sigma^r_x(R)$ to $\sigma^r_y(R')$ for each $0 \leq r \leq p$. We
say that $R$ is {\bf locally equivalent} to $R'$ if there exists
an open subset $U \subset M$ (resp. $ U' \subset M'$) and a line
bundle isomorphism $\varphi: E|_U \rightarrow E'|_{U'}$ such that
the induced bundle isomorphism $J^p(\varphi): J^p(E)|_{U}
\rightarrow J^p(E')|_{U'}$ maps $R|_{U}$ onto $R'|_{U'}$.

\medskip Many interesting examples of  integrable systems of linear
differential equations of finite order arises from the projective
differential geometry (e.g. [W]). In modern terms, it can be
presented as follows.  Let $M \subset {\bf P}_N$ be a
non-degenerate complex sub-manifold. Let $E$ be the line bundle on
$M$ defined by the restriction of the hyperplane line bundle of
${\bf P}_N$ and $\Sigma$ be the $(N+1)$-dimensional space of
sections of $E$ coming from the sections of the hyperplane line
bundle. Let $p>0$ be a positive integer such that an element of
$\Sigma$ having a zero of order $\geq p$ at a point of $M$ is
identically zero. Define $R_x \subset J^p_x(E)$ by
$$ R_x := \{ j^p_x(s) | s \in \Sigma \}.$$ Then over a
Zariski-open subset $M^o$ of $M$, $$R := \bigcup_{x \in M^o} R_x
$$ defines a sub-bundle of $J^p(E)|_{M^o}$, called the {\bf
hyperplane system} on $M$. By the construction, $R$ is of finite
type and integrable. In fact, $$\pi^p_q (R_x) = \{ j^q_x(s) |s \in
\Sigma \}.$$

\medskip
{\bf Proposition 1} {\it  Let $M$ and $ M'$ be two non-degenerate
sub-manifolds of ${\bf P}_N$.  Let $R$ (resp. $ R'$) be the
hyperplane system on $M$ (resp. $M'$). Assume that $R$ and $R'$
are of the same order. Then $R$ and $R'$ are locally equivalent if
and only if there exists a projective linear transformation
sending an open subset of $M$ onto an open subset of $M'$.}

\medskip
{\it Proof}. The `if' part is trivial. Let us prove the `only if'
part.  Let $E$ (resp. $E'$) be the restriction of the hyperplane
line bundle on $M$ (resp. $M'$). Let $\Sigma$ (resp. $\Sigma'$) be
the sections of $E$ (resp. $E'$) coming from the global sections
of the hyperplane line bundle on ${\bf P}_N$. If $R$ and $R'$ are
locally equivalent, there exists an open subset $U \subset M$
(resp. $U' \subset M'$) and an isomorphism of line bundles
$\varphi: E|_U \rightarrow E'|_{U'}$ mapping $R|_U$ onto
$R'|_{U'}$. The solutions of $R|_U$ are sent to the solutions of
$R'|_{U'}$. By the assumption that $R$ and $R'$ are of finite
type,   every local solution of $R$ (resp. $R'$) belongs to
$\Sigma$ (resp. $\Sigma'$). Thus $\varphi$ induces a linear
isomorphism from $\Sigma$ to $\Sigma'$ which gives a projective
linear transformation of ${\bf P}_N$ sending $U$ to $U'$. $\Box$

\medskip
For $x \in M^o$, the symbol $\sigma^r_x(R) \subset {\rm Sym}^r
T_x^*$ is called the {\bf r-th fundamental form} of $M$ at $x$ and
is denoted by $F^r_x(M)$. We decree $F^r_x(M) =0$ for $r \geq p$.
For two sub-manifolds $M, M' \subset {\bf P}_N$, we say that the
fundamental forms of $M$ at $x \in M$ are isomorphic to the
fundamental forms of $M'$ at $y\in M'$, if there exists a linear
isomorphism $\theta: T_x \rightarrow T_y$ such that the induced
linear map $\theta^r: {\rm Sym}^r T^*_x \rightarrow {\rm Sym}^r
T^*_y$ sends $F^r_x(M)$ to $F^r_y(M')$ for each $r \geq 0$.

\medskip
{\bf Remark} Usually, the $r$-th fundamental form of $M$ at $x$ is
defined as a linear map ${\rm Sym}^r T_x \rightarrow N^r_x$ where
$N^r_x$ is the $r$-th normal space at $x$ (e.g. [GH]). Our
$F^r_x(M) \subset {\rm Sym}^r T^*_x$ is the image of the dual of
this linear map.

\section{Fundamental forms of Hermitian symmetric spaces}
We start with a few definitions regarding  representations of
abelian Lie algebras on graded vector spaces.

Let $W = \bigoplus_{r=0}^p W_r$ be a finite-dimensional graded
vector space with $\dim W_0=1$. For each integer $k$, define
$${\bf gl}(W)_k := \{ X \in {\bf gl}(W) \; |\; X(W_r) \subset
W_{k+r} \mbox{ for all } r\}. $$ Let $V$ be a finite dimensional
vector space and $\mu:V \rightarrow {\bf gl}(W)$ be a
representation of the abelian Lie algebra $V$ such that $\mu(V)
\subset {\bf gl}(W)_{-1}$. Setting ${\bf p}_{-1} = \mu(V)$, define
inductively for each $k \geq 0$,
$${\bf p}_k := \{ X \in {\bf gl}(W)_k \; | \; [{\bf p}_{-1}, X]
\subset {\bf p}_{k-1} \}.$$ Then the graded Lie subalgebra ${\bf
p}= \bigoplus_{k=-1}^{\infty} {\bf p}_k$ of ${\bf gl}(W)$ is
called the {\bf prolongation} of $\mu$. For each non-negative
integer $r$, define a homomorphism $\chi^{\mu}_r: W_r \rightarrow
{\rm Sym}^r V^* \otimes W_0$ by
$$\chi^{\mu}_r(w) (X_1, \ldots, X_r) :=
 \mu(X_1) \cdots \mu(X_r) (w)$$ for $w \in W_r$ and $X_1,
\ldots, X_r \in V$. The right hand side is symmetric in $X_1,
\ldots, X_r$  because $\mu$ is a representation of the abelian Lie
algebra $V$. After fixing an identification $W_0 = {\bf C}$, we
denote by $\chi^{\mu}: W \rightarrow \bigoplus_{r=0}^p {\rm Sym}^r
V^*$ the graded vector space homomorphism obtained as the sum of
$\bigoplus_{r=0}^p \chi^{\mu}_r$. The homomorphism $\chi^{\mu}$ is
called the {\bf symbol map} of the representation $\mu: V
\rightarrow {\bf gl}(W)_{-1}$. The following proposition is
obvious.

\medskip
{\bf Proposition 2} {\it For $X \in V$, let $\iota(X): {\rm
Sym}^{\bullet}V^* \rightarrow {\rm Sym}^{\bullet} V^*$ denote the
inner multiplication by $X$. Then the image of $\chi^{\mu}$ is
invariant under $\iota(X)$ for any $X \in V$ and $$\iota(X) \circ
\chi^{\mu}(w) = \chi^{\mu}(\mu(X) \cdot w) $$ for every $X \in V$
and $w \in W$.}

\medskip
Now  let ${\bf l} = {\bf l}_{-1} + {\bf l}_0 + {\bf l}_1$ be a
semi-simple graded Lie algebra. Let $L$ be a connected Lie group
with Lie algebra ${\bf l}$ and $L'$ be a connected subgroup of $L$
with Lie algebra ${\bf l}' = {\bf l}_0 + {\bf l}_1$. Then the
homogeneous space $L/L'$ is a compact Hermitian symmetric space.
Conversely, for any compact Hermitian symmetric space, let ${\bf
l}$ be the Lie algebra of global holomorphic vector fields and
${\bf l}'$ be the subalgebra consisting of vector fields vanishing
at a base point. Then there exists a grading  ${\bf l}= {\bf
l}_{-1} + {\bf l}_0 + {\bf l}_1$ such that ${\bf l}'= {\bf l}_0
\oplus {\bf l}_1$.

\medskip
 Given a non-degenerate
equivariant embedding of the Hermitian symmetric space $L/L'
\subset {\bf P}S^*$ for a finite-dimensional vector space $S$, we
get a faithful irreducible representation $\rho: {\bf l}
\rightarrow {\bf gl}(S)$ such that
$$S_0 := \{ s \in S| \rho({\bf l}_{-1}) (s) =0 \}$$ is
1-dimensional. Define for $r \geq 0$,
$$S_{r+1} := \rho({\bf l}_1) (S_r). $$ There exists $p>0$ such
that $S_r \neq 0$ for $r=0, 1, \ldots, p-1$ and $S_r =0$ for $r
\geq p$. Moreover $$S = \bigoplus_{r=0}^p S_r.$$ Set $V:= {\bf
l}_{-1}$ and let $ \nu: V \rightarrow {\bf gl}(S)$ be the
representation obtained by restricting $\rho:{\bf l} \rightarrow
{\bf gl}(S) $ to ${\bf l}_{-1}$.  Then $\nu: V \rightarrow {\bf
gl}(S)_{-1}$ with respect to the grading $S = \bigoplus_{r=0}^p
S_r$.

\medskip
We recall the following result on the prolongation of $\nu$.

\medskip
{\bf Proposition 3 ([S,Proposition 4.4.1], or [SYY, (2.5)] for the
irreducible case)} {\it Regard ${\bf l}$ as a subalgebra of ${\bf
gl}(S)$ by the given faithful irreducible representation
$\rho:{\bf l} \rightarrow {\bf gl}(S)$. Let ${\bf z}$ be the
centralizer  of ${\bf l}$ in  ${\bf gl}(S)$.   For the
representation $\nu: V := {\bf l}_{-1} \rightarrow {\bf
gl}(S)_{-1}$, the prolongation is ${\bf p}_{-1} \cong V$ (i.e.
$\nu$ is injective), ${\bf p}_0 = {\bf l}_0 + {\bf z}$, ${\bf p}_1
= {\bf l}_1$ and ${\bf p}_k =0$ for $k \geq 2$. In particular,
${\bf p}$ is a reductive graded Lie algebra whose semi-simple part
is the graded semi-simple Lie algebra ${\bf l}$. }

\medskip
The following result on the symbol map of $\nu$ is contained in
[S].  Note that ${\bf l}_{-1}$ is naturally identified with the
tangent space of the Hermitian symmetric space $L/L'$ at the base
point.

\medskip
{\bf Proposition 4 } {\it Let $V ={\bf l}_{-1}$ and $\nu:V
\rightarrow {\bf gl}(S)_{-1}$ be the representation obtained by
restricting $\rho:{\bf l} \rightarrow {\bf gl}(S)$ with respect to
the grading $S = \bigoplus_{r=0}^p S_r$. Then $\chi^{\nu}:
\bigoplus_{r=0}^p S_r \rightarrow \bigoplus_{r=0}^p {\rm Sym}^r
V^*$ is injective and the image $\chi^{\rho}_r(S_r) \subset {\rm
Sym}^r V^*$ coincides with the $r$-th fundamental form
$F^r_y(L/L') \subset {\rm Sym}^r ({\bf l}_{-1})^*$ of $L/L'
\subset {\bf P}S^*$. }

\medskip
{\it Proof}. The injectivity of $\chi_r^{\nu}$ is [S, Proposition
4.3.1]. That the image $\chi^{\rho}_r(S_r) \subset {\rm Sym}^r
V^*$ coincides with the $r$-th fundamental form is [S, Proposition
2.4.1.]. $\Box$

\medskip
Now we regard $S=\bigoplus_{r=0}^pS_r$ as a graded vector subspace
of $\bigoplus_{r=0}^p {\rm Sym}^r V^*$ by the injection
$\chi^{\nu}$. By Proposition 2, the action of $V$ on $S$ by $\nu$
coincide with the restriction of the inner multiplication. The
following invariant-theoretic result implies the `rigidity' of the
fundamental forms of Hermitian symmetric spaces.

\medskip
{\bf Proposition 5} {\it Regard $S= \bigoplus_{r=0}^p S_r$ as a
graded subspace of $ \bigoplus_{r=0}^p {\rm Sym}^r V^*$ as
explained above.  Let $P = \bigoplus_{r=0}^p P_r$ be a graded
subspace of $ \bigoplus_{r=0}^p {\rm Sym}^r V^*$ which has the
same discrete invariants as $S$ under the action of $GL(V)$. Then
$P$ is in the $GL(V)$-orbit of  $S$ .  In other words,
 there exist finitely many $GL(V)$-invariant
${\bf Z}$-valued functions on the  Grassmannian  $Gr(\dim S,
\bigoplus_{r=0}^p {\rm Sym}^r V^*)$  of $(\dim S)$-dimensional
subspaces in $\bigoplus_{r=0}^p{\rm Sym}^r V^*$, such that if all
these functions have the same values on $P$ and $S$, then $P$ is
in the orbit of $S$ under the $GL(V)$-action on the Grassmannian.}

\medskip
{\it Proof}. By Proposition 2, $S$ is invariant under $\iota(X)$
for all $X \in V$. Since $\dim \iota(X) P$ and $\dim (P \cap
\iota(X) P)$ must be equal to the corresponding dimensions for
$S$, we conclude that $P$ is also invariant under $\iota(X)$.
Moreover the action of $\iota(X)$ and $\iota(Y)$ on $P$ commute
for any $X, Y \in V$. Thus we have a representation $\mu: V
\rightarrow {\bf gl}(P)_{-1}$. The dimension of the prolongation
of $\mu$ and the dimension of the prolongation of $\nu:V
\rightarrow {\bf gl}(S)$ are equal by the assumption on $P$. So
are the dimensions of their radicals. This implies, by Proposition
3, the prolongation ${\bf p}= {\bf p}_{-1} + {\bf p}_0 + {\bf p}_1
+ \cdots$ of $\mu$ is reductive and isomorphic to ${\bf l}_{-1} +
( {\bf l}_0 + {\bf z}) + {\bf l}_1$. By the same reasoning, the
action of the prolongation ${\bf p}$ on $P$ is isomorphic to the
irreducible representation of ${\bf l}$ on $S$. In particular,
$\mu: V \rightarrow {\bf gl}(P)_{-1}$ is isomorphic to $\nu: V
\rightarrow {\bf gl}(S)_{-1}$  under a linear automorphism of $V$
and a graded vector space isomorphism $P \cong S$. Consider
$\chi^{\mu}: P \rightarrow \bigoplus_{r=0}^p {\rm Sym}^r V^*$.
 The  image $\chi^{\mu}(P) \subset \bigoplus_{r=0}^p{\rm Sym}^r
V^*$ is completely determined by the structure of $P$ as a
$V$-module. Thus $\dim \chi^{\mu}(P) = \dim \chi^{\nu}(S)$ and
$\chi^{\mu}(P)$ is in the $GL(V)$-orbit of $\chi^{\nu}(S)$.  But
by our definition of $S$ as the subspace of $\bigoplus_{r=0}^p
{\rm Sym}^r V^*$, $\chi^{\nu} (S) =S$. Since $\dim \chi^{\mu}(P)
\cap P$ is an integer-valued invariant of $P$, we conclude $P =
\chi^{\mu}(P)$ is in the $GL(V)$-orbit of  $ \chi^{\nu} (S) = S$.
$\Box$

\section{Vanishing of Spencer cohomology}

Let ${\bf l}={\bf l}_{-1} + {\bf l}_0 + {\bf l}_1$ be a graded
semi-simple Lie algebra associated to a Hermitian symmetric space
as explained in Section 3. We recall the description of the
gradation of ${\bf l}$ in terms of the root system. Let ${\bf h}
\subset {\bf l}_0$ be a Cartan subalgebra of ${\bf l}$. Then there
exists a unique element $Z \in {\bf h}$, called the {\bf
characteristic element}, with the property that
$${\bf l}_i = \{ X \in {\bf l} \;|\; [Z, X] = i\; X\}$$ for each $i =
-1, 0, 1$. Choose a simple root system $\Delta = \{ \alpha_1,
\ldots, \alpha_{\ell} \}$ such that $\alpha_i(Z) \geq 0$ for all
$i = 1, \ldots, \ell$. Then $Z$ determines a partition $\Phi^{+} =
\Phi^+_0 \cup \Phi^+_1$ of the set of positive roots by $\Phi^+_k
= \{\alpha \in \Phi^+| \alpha(Z) =k \}$ for $k=0,1,$ such that
\begin{eqnarray*} {\bf l}_{-1} &=& \bigoplus_{\alpha \in \Phi^+_1}
{\bf g}_{-\alpha}\\
 {\bf l}_0 &=& {\bf h} \oplus
\bigoplus_{\alpha \in \Phi^+_0} ({\bf g}_{\alpha} \oplus {\bf
g}_{-\alpha}) \\
 {\bf l}_1 &=& \bigoplus_{\alpha \in \Phi^{+}_1}
{\bf g}_{\alpha} \end{eqnarray*} where ${\bf g}_{\alpha}$ is the
root space. Define $$\Delta_1 = \{ \alpha_{i_1}, \ldots,
\alpha_{i_m} \} := \{ \alpha \in \Delta| \alpha(Z) =1 \}.$$ Recall
$$<\!\beta, \gamma\!> = \frac{2 (\beta, \gamma)}{(\gamma, \gamma)}$$
for $\beta, \gamma \in {\bf h}^*$. Let $\lambda_1, \ldots,
\lambda_{\ell}$ be the fundamental weights and $h_{\alpha_1},
\ldots, h_{\alpha_{\ell}}$ be the co-roots, namely,
$$<\!\lambda_i, \alpha_j \!> \; = \; \delta_{ij}, \;\;\;
\beta(h_{\alpha_i}) =  \; <\! \beta, \alpha_i \!> \;\;\; \mbox{
for all } \beta \in {\bf h}^*.$$ By a direct computation, we have

\medskip
{\bf Lemma 1} {\it Suppose $\lambda_{i_1} + \cdots + \lambda_{i_m}
= \sum_{i=1}^{\ell} a_i \alpha_i.$ Then $Z = \sum_{i=1}^{\ell} a_i
h_{\alpha_i}. $}

\medskip
Let $W$ be the Weyl group of the root system $\Phi$. For an
element $\sigma \in W$, put  $\Phi_{\sigma} = \sigma(-\Phi^+) \cap
\Phi^+$. Define
\begin{eqnarray*} W^0 &:=& \{
\sigma \in W \;| \; \Phi_{\sigma} \subset \Phi^+_1 \} \\
W^0(1) &:=& \{ \sigma \in W^0 \; | \; \sharp(\Phi_{\sigma}) =1 \}
\end{eqnarray*} where $\sharp(\Phi_{\sigma})$ is the number of
roots in $\Phi_{\sigma}$. For $\alpha_{i} \in \Delta$, let
$\sigma_i \in W$ be the corresponding reflection.

\medskip
{\bf Lemma 2}  $W^0(1) = \{ \sigma_{\alpha} \; | \; \alpha \in
\Delta_1 \} = \{ \sigma_{i_1}, \ldots, \sigma_{i_m} \}.$

\medskip
{\it Proof}. The inclusion $\{ \sigma_{i_1}, \ldots, \sigma_{i_m}
\} \subset W^0(1)$ is immediate. By [H, 10.3 Lemma A], elements of
$W^0(1)$ are simple reflections. If $\sigma_i \in W^0(1)$, then
$\Phi_{\sigma_i} = \{ \sigma_i \}$ and $\Phi_\sigma \subset
\Phi^+_1$. Thus $\alpha_i \in \Delta_1 $. $\Box$

\medskip
For a finite-dimensional ${\bf l}$-module $\Gamma$ for which the
weights of $Z$ are integers, let $\Gamma_p$ be the subspace of
$\Gamma$ consisting of eigenvectors of $Z$ with eigenvalue $p \in
{\bf Z}$. Set $C^{p,q}  := \Lambda^q ({\bf l}_{-1})^* \otimes
\Gamma_{p-1}$ and define $\partial: C^{p,1} \rightarrow C^{p-1,
2}$ by
$$\partial c (X_0, X_1) = X_0 \cdot c(X_1) - X_1 \cdot c(X_0) $$
for $c \in C^{p,1} = ({\bf l}_{-1})^*\otimes \Gamma_{p-1}$ and
$\partial: C^{p+1, 0} \rightarrow C^{p,1}$ by $$\partial c(X) = X
\cdot c$$ for $x \in \Gamma_p$. For each integer $p$, the $p$-th
{\bf Spencer cohomology group} associated to $\Gamma$ is defined
by
$$H^{p,1}({\bf l}_{-1}, \Gamma) : = \frac{({\rm Ker}
\partial: C^{p,1} \longrightarrow  C^{p-1,2})}{({\rm Im} \partial:
C^{p+1,0} \longrightarrow C^{p,1})}.$$ They give the decomposition
of the Lie algebra cohomology group
$$H^1({\bf l}_{-1}, \Gamma) = \oplus_p H^{p,1}({\bf l}_{-1}, \Gamma)$$
which is just the $Z$-eigenspace decomposition.

\medskip
Given a faithful irreducible ${\bf l}$-module $S$ with $\dim
S_0=1$, ${\bf g} \subset {\bf gl}(S)$ be the prolongation of ${\bf
l}_{-1}$-action on $S$ as defined in the last section. Let ${\bf
g}^{\perp}\subset {\bf gl}(S)$ be the orthogonal complement with
respect to the Killing form on ${\bf gl}(S)$. Then ${\bf
g}^{\perp}$ is an ${\bf l}$-module where the weights of $Z$ are
integers. So $H^{p,1}({\bf l}_{-1}, {\bf g}^{\perp})$ makes sense.

\medskip
  Let $R$ be a system of linear differential equations
of order $p$ on a complex manifold $M$. Given a graded subspace
$S= \bigoplus_{i=0}^p S_i \subset {\rm Sym}^{\bullet} V$ for a
vector space $V$, $\dim V = \dim M$, we say that the symbol
$\sigma_x(R)$ is isomorphic to $S$ if there exists a linear
isomorphism $\psi: T^*_x \rightarrow V$ such that the induced
graded linear isomorphism ${\rm Sym}^{\bullet}T_x^* \rightarrow
{\rm Sym}^{\bullet} V$ sends $\sigma_x(R)$ onto $S$.

We will reduce the proof of  Theorem 1 and Theorem 2 to the
vanishing of certain Spencer cohomology groups by using the
following result of Y. Se-ashi.  Theorem 3 appeared as Theorem
5.1.2, Theorem 5.2.2 and Theorem 5.3.1 in [S]. We recommend [SYY,
2.4] for an over-view of its proof, where it appeared as `Theorem
A'.

\medskip
{\bf Theorem 3 (Se-ashi) } {\it Let $H \subset {\bf P}_N$ be an
equivariantly embedded Hermitian symmetric space and $F =
\bigoplus_{r=0}^{p} F^r $ be the fundamental forms at a point of
$H$. Assume that $H^{p,1}({\bf l}_{-1}, {\bf g}^{\perp}) = 0$ for
all $p \geq 1$.  Suppose a system $R$ on a line bundle $E$ of a
complex manifold $M$ has its symbol at each point $x \in M$
isomorphic to $F$. Then $R$ is locally equivalent to the
hyperplane system on $H \subset {\bf P}_N.$}

\medskip
For the computation of the Spencer cohomology we will use
Kostant's result.  The following is a translation of [K, Theorem
5.14] into our situation (cf. [SYY,2.5 Theorem B]). See also
[Y,5.1] for a stream-lined review of Kostant theory. Note that
${\bf l}_1= ({\bf l}_{-1})^*$ by the Killing form.

\medskip
{\bf Theorem 4 (Kostant) } {\it Let $\Gamma$ be an irreducible
${\bf l}$-module and  $\gamma$ be its lowest weight.   Then we
have an ${\bf l}_0$-module decomposition
$$H^{1}({\bf l}_{-1}, \Gamma) = \bigoplus_{\sigma \in W^0(1)} {\cal H}^{\xi_{\sigma}}$$
where ${\cal H}^{\xi_{\sigma}}$ is the irreducible ${\bf
l}_0$-module with the lowest weight $\xi_{\sigma} :=
\sigma(\gamma) + \alpha$ for  the simple root $\alpha$ satisfying
$ \sigma = \sigma_{\alpha}$. }

\medskip
This has the following consequence for the Spencer cohomology
groups.

\medskip
{\bf Proposition 6} {\it Let  $Z \in {\bf h}$ be the
characteristic element of the graded Lie algebra ${\bf l}$. Let
$\Gamma$ be an irreducible ${\bf l}$-module such that the weights
of $Z$ on $\Gamma$ are integers. Let $\gamma$ be the lowest weight
of $\Gamma$. Then $H^{p, 1}({\bf l}_{-1}, \Gamma) \neq 0$ only if
$p-1 = \sigma_{i_k}(\gamma)(Z)$ for some $i_k$. }

\medskip
{\it Proof}. $H^{p,1}({\bf l}_{-1}, \Gamma)$ is the eigenspace of
$Z$ with eigenvalue $p$. By Theorem 4, it is non-zero if there
exists $\sigma_i \in W^0(1)$ such that $$p= \xi_{\sigma_i}(Z) =
\sigma_i(\gamma)(Z) + \alpha_i(Z). $$ Thus the result follows from
$W^0(1) = \{ \sigma_{i_1}, \ldots, \sigma_{i_m} \}$ by Lemma 2 and
$\alpha_{i_k}(Z) =1$. $\Box$

\medskip
The following is an honest generalization of [SYY, Theorem 2].

\medskip
{\bf Proposition 7} {\it  Let ${\bf l}= {\bf l}_{-1} + {\bf l}_0 +
{\bf l}_1$ be the graded semi-simple Lie algebra associated to a
Hermitian symmetric space which contains neither a projective
space nor a hyperquadric as an irreducible factor.  Then for any
irreducible ${\bf l}$-module $\Gamma$,
 $$H^{p,
1}({\bf l}_{-1}, \Gamma) = 0 \mbox{ for all } p \geq 1.$$}

{\it Proof}. \begin{eqnarray*} \sigma_{\i_k}(\gamma)(Z) &=&
(\gamma -<\!\gamma, \alpha_{i_k}\!> \alpha_{i_k})(Z) \\
&=& \gamma(Z) -<\!\gamma, \alpha_{i_k}\!>. \end{eqnarray*}  By
Lemma 1,
\begin{eqnarray*} \sigma_{i_k}(\gamma)(Z) &=& (a_{i_k} -1) <\!\gamma,
\alpha_{i_k}\!> + \sum_{i \neq i_k} a_i <\! \gamma, \alpha_i \!>.
\end{eqnarray*}
Since $\gamma$ is the lowest weight, we have $<\!\gamma,
\alpha_i\!> \; \leq 0$ for $i=1, \ldots, \ell$ and  $<\!\gamma,
\alpha_j\!> \;<0$ for some $j$. Since $a_i$ is an entry of the
inverse matrix of the Cartan matrix, $a_i
>0$. In particular, if $a_{i_k}
>1$, then $\sigma_{i_k}(\gamma)(Z) <0$ for every $\gamma$. The
only cases when $a_{i_k} \leq 1$ are the following (e.g. [H, 13.2
Table 1]).
\begin{itemize}
\item{$a_{i_k} = \frac{\ell}{\ell +1}$ when $\alpha_{i_k}$ is the root $\alpha_1$ of $A_{\ell}, \ell \geq 1$. }
\item{$a_{i_k} =1$ when $\alpha_{i_k}$ is the root $\alpha_2$ of $A_3$.}
\item{$a_{i_k} =1$ when $\alpha_{i_k}$ is the root $\alpha_1$ of $B_{\ell}, \ell \geq 2$.}
\item{ $a_{i_k} =1$ when $\alpha_{i_k}$ is the root $\alpha_1$ of $D_{\ell},
\ell \geq 4$.} \end{itemize} In these cases, the corresponding
irreducible factor of the Hermitian symmetric space is a
projective space or a hyperquadric. $\Box$

\medskip
{\it Proof of Theorem 1}. Applying Proposition 7 to each
irreducible factor of ${\bf g}^{\perp}$, we see that $H^{p,1}({\bf
l}_{-1}, {\bf g}^{\perp}) =0$ for all $p \geq 1$.  By Proposition
5, the fundamental forms of $M$ at each point of a neighborhood of
$x$ are isomorphic to the fundamental forms of $H$ at a point.
Thus Theorem 1 follows from  Theorem 3 and Proposition 1. $\Box$

\medskip
To prove Theorem 2, let us study the corresponding ${\bf
l}$-modules.  Let ${\bf l} = {\bf sl}_{n_1+1} \oplus \cdots \oplus
{\bf sl}_{n_m+1}$ with $n_i \geq 2$ for all $i=1, \ldots, m$. Let
us denote the simple roots of ${\bf sl}_{n_j+1}$ by $\{
\alpha^j_1, \ldots, \alpha^j_{n_j} \}$ using the standard
enumeration of simple roots of the simple Lie algebra of type $A$.
For the Hermitian symmetric space ${\bf P}_{n_1} \times \cdots
\times {\bf P}_{n_m}$, the corresponding gradation of ${\bf l}$
has
$$\Delta_1 = \{ \alpha^1_1, \alpha^2_1, \ldots, \alpha^m_1 \}.$$
The fundamental weights of ${\bf l}$ are denoted by $$\lambda^1_1,
\ldots, \lambda^1_{n_1}, \lambda^2_1, \ldots, \lambda^2_{n_2},
\ldots, \lambda^m_1, \ldots, \lambda^m_{n_m}.$$ Let $U^j = {\bf
C}^{n_j+1}$ be the basic representation space of ${\bf
sl}_{n_j+1}$ corresponding to the fundamental weight
$\lambda^j_1$.  Let $S:=(U^1)^* \otimes \cdots \otimes (U^m)^*$ be
the basic representation space of ${\bf l}$. Using ${\bf
gl}((U^i)^*) = {\bf C} \cdot {\rm Id} + {\bf sl}((U^i)^*)$ and
${\bf sl}((U^i)^*)\cong{\bf sl}(U^i)$, we have the ${\bf
l}$-module decomposition
\begin{eqnarray*} {\bf gl}(S) &=&
\left( {\bf C} \cdot {\rm Id} +{\bf sl}(U^1) \oplus \cdots \oplus
{\bf sl}(U^m) \right) + {\bf g}^{\perp}
\end{eqnarray*} where ${\bf g}^{\perp}$ is a sum of  tensor products of
the form ${\bf sl}(U^{i_1}) \otimes \cdots \otimes {\bf
sl}(U^{i_k})$ for some $k \geq 2$. By Proposition 1, Proposition 5
and Theorem 3, the proof of Theorem 2 is reduced to the following.

 \medskip
 {\bf Proposition 8} {\it In the  above notation,
 $H^{p,1}({\bf l}_{-1}, {\bf g}^{\perp}) =0 \mbox{ for all } p \geq 1$.}

\medskip {\it Proof}.
  Let $\gamma$ be the
lowest weight of ${\bf sl}(U^{i_1}) \otimes \cdots \otimes {\bf
sl}(U^{i_k}), k \geq 2$. By Proposition 6, it suffices to show,
for each $\ell=1, \ldots, m$,
$$\gamma(Z) -<\!\gamma, \alpha^{\ell}_1\!> \; \leq -1.$$ Since the
 adjoint representation ${\bf sl}(U^j)$ has highest weight
$\lambda^j_1 + \lambda^j_{n_j}$,
$$\gamma = - \sum_{j=1}^k (\lambda^{i_j}_1 + \lambda^{i_j}_{n_{i_j}}).$$
When $\lambda^j_1 = \sum _{i=1}^{n_j} a^j_i \alpha^j_i$, Lemma 1
says
$$Z = \sum_{j=1}^m \sum_{i=1}^{n_j} a^j_i h_{\alpha^j_i} .$$ Using
 (e.g. from [H, 13.2 Table 1])  $a^j_1 = \frac{n_{j}}{n_j+1}$ and
$a^j_{n_j} = \frac{1}{n_j+1}$, we have
$$\gamma (Z) = -\sum_{j=1}^k (a^{i_j}_1 + a^{i_j}_{n_{i_j}}) = -k$$ while
$$ -<\!\gamma, \alpha^{\ell}_1 \!> \;=\; <\! \lambda^{i_j}_1 + \lambda^{i_j}_{n_{i_j}},
\alpha^{i_j}_1 \!> \;= \; 1 \mbox{ when } \ell = i_j \mbox{ and
}$$ $$- <\! \gamma, \alpha^{\ell}_1 \!> \;=0 \mbox{ when } \ell
\neq i_1, \ldots, i_k.$$  Thus $\gamma(Z) -<\!\gamma,
\alpha^{\ell}_1\!> \; \leq -1$. $\Box$

\medskip
{\bf Remark} When $n_1 =1$ and $2 \leq n_2, \ldots, n_m$, the
adjoint representation ${\bf sl}(U^1)$ has highest weight $2
\lambda^1_1 = \alpha^1_1$. Thus ${\bf g}^{\perp}$ has an
irreducible factor with lowest weight
$$\gamma = - 2 \lambda_1 - \lambda^2_1 -
\lambda^2_{n_2}.$$ Then the same calculation as in the proof of
Proposition 8 shows
$$\gamma (Z) = -2 a^1_1 -a^2_1 -a^2_{n_2} = -2$$
$$-<\! \gamma, \alpha^1_1\!> = 2.$$ Thus $H^{1,1}({\bf l}_{-1},
{\bf g}^{\perp}) \neq 0$.  It is not known whether ${\bf P}_1
\times {\bf P}_{n_2} \times \cdots {\bf P}_{n_m}$ with $2 \leq
n_2, \ldots n_m$  is determined by its fundamental forms. For
$m=2$, this was conjectured to be true in [L1, Section 6].

\bigskip
{\bf References}

\medskip

[GH] Griffiths, P. and Harris, J.: Algebraic geometry and local
differential geometry. Ann. scient. \'Ec. Norm. Sup. {\bf 12}
(1979) 355-432

[H] Humphreys, J.: {\it Introduction to Lie algebras and
representation theory}. Grad. Texts Math. {\bf 9}, Springer, New
York-Heidelberg-Berlin 1972

 [JM] Jensen, G. and Musso, E.: Rigidity of hypersurfaces in
complex projective space. Ann. scient. \'Ec. Norm. Sup. {\bf 27}
(1994) 227-248

 [K] Kostant, B.: Lie algebra cohomology and generalized
Borel-Weil theorem. Ann. Math. {\bf 74} (1961) 329-387

[L1] Landsberg, J. M.: On the infinitesimal rigidity of
homogeneous varieties. Compositio Math. {\bf 118} (1999) 189-201

[L2] Landsberg, J. M.: Griffiths-Harris rigidity of compact
Hermitian symmetric spaces. preprint. alg-geom/0207287.

 [SYY] Sasaki, T., Yamaguchi, K. and
Yoshida, M.: On the rigidity of differential systems modelled on
Hermitian symmetric spaces and disproofs of a conjecture
concerning modular interpretations of configuration spaces. Adv.
Stud.  Pure Math. {\bf 25} {\it CR-geometry and overdetermined
systems.} (1997) 318-354

[S] Se-ashi, Y., On differential invariants of integrable finite
type linear differential equations. Hokkaido Math. J. {\bf 17}
(1988) 151-195

[W] Wilczynski, E. J.: {\it Projective differential geometry of
curves and ruled surfaces}. reprinted by Chelsea Publ. Co.
Teubner, 1906

 [Y] Yamaguchi, K.: Differential systems associated with
simple graded Lie algebras. Adv. Stud. Pure Math. {\bf 22} {\it
Progress in differential geometry}. (1993) 413-494

\bigskip
\begin{tabbing}
aaaaaaaaaaaaaaaaaaaaaaaaaaaaaaaaaaaaaaaa\= \kill
Jun-Muk Hwang \> Keizo Yamaguchi \\

School of Mathematics \> Department of Mathematics\\

Korea Institute for Advanced Study  \>  Graduate School of Science \\

207-43 Cheongnyangni-dong \> Hokkaido University\\

Seoul 130-722, Korea  \> Sapporo 060-0810, Japan\\

e-mail: jmhwang@ns.kias.re.kr \> e-mail:
yamaguch@math.sci.hokudai.ac.jp
\end{tabbing}

\end{document}